\documentstyle[12pt]{article}

\begin{document}

\newtheorem{theorem}{Theorem}[section]
\newtheorem{corollary}[theorem]{Corollary}
\newtheorem{definition}[theorem]{Definition}
\newtheorem{conjecture}[theorem]{Conjecture}
\newtheorem{question}[theorem]{Question}
\newtheorem{lemma}[theorem]{Lemma}
\newtheorem{proposition}[theorem]{Proposition}
\newtheorem{example}[theorem]{Example}
\newtheorem{problem}[theorem]{Problem}
\newenvironment{proof}{\noindent {\bf
Proof.}}{\rule{3mm}{3mm}\par\medskip}
\newcommand{\remark}{\medskip\par\noindent {\bf Remark.~~}}
\newcommand{\pp}{{\it p.}}
\newcommand{\de}{\em}

\newcommand{\JEC}{{\it Europ. J. Combinatorics},  }
\newcommand{\JCTB}{{\it J. Combin. Theory Ser. B.}, }
\newcommand{\JCT}{{\it J. Combin. Theory}, }
\newcommand{\JGT}{{\it J. Graph Theory}, }
\newcommand{\ComHung}{{\it Combinatorica}, }
\newcommand{\DM}{{\it Discrete Math.}, }
\newcommand{\ARS}{{\it Ars Combin.}, }
\newcommand{\SIAMDM}{{\it SIAM J. Discrete Math.}, }
\newcommand{\SIAMADM}{{\it SIAM J. Algebraic Discrete Methods}, }
\newcommand{\SIAMC}{{\it SIAM J. Comput.}, }
\newcommand{\ConAMS}{{\it Contemp. Math. AMS}, }
\newcommand{\TransAMS}{{\it Trans. Amer. Math. Soc.}, }
\newcommand{\AnDM}{{\it Ann. Discrete Math.}, }
\newcommand{\NBS}{{\it J. Res. Nat. Bur. Standards} {\rm B}, }
\newcommand{\ConNum}{{\it Congr. Numer.}, }
\newcommand{\CJM}{{\it Canad. J. Math.}, }
\newcommand{\JLMS}{{\it J. London Math. Soc.}, }
\newcommand{\PLMS}{{\it Proc. London Math. Soc.}, }
\newcommand{\PAMS}{{\it Proc. Amer. Math. Soc.}, }
\newcommand{\JCMCC}{{\it J. Combin. Math. Combin. Comput.}, }
\newcommand{\GC}{{\it Graphs Combin.}, }

\title{The Maximum Wiener Index of Trees with Given  Degree Sequences\thanks{
 Supported by National Natural Science Foundation of China
(No.10531070), National Basic Research Program of China 973 Program
(No.2006CB805900), National Research Program of China 863 Program
(No.2006AA11Z209) and a grant of Science and Technology  Commission
of Shanghai Municipality (STCSM No: 09XD1402500).\newline \indent
$^{\dagger}$Correspondent author: Xiao-Dong Zhang (Email:
xiaodong@sjtu.edu.cn)}}
\author{  Xiao-Dong Zhang$^{\dagger}$, Yong Liu  and Min-Xian Han  \\
{\small Department of Mathematics}\\
{\small Shanghai Jiao Tong University} \\
{\small  800 Dongchuan road, Shanghai, 200240, P.R. China}\\
 }
\date{}
\maketitle
 \begin{abstract}
 The Wiener index  of a  connected graph is
 the sum of topological distances between all pairs of vertices.
 Since Wang in \cite{wang2008} gave a mistake result on the maximum Wiener index
 for given  tree
 degree sequence, in this paper, we investigate the maximum Wiener index of  trees with given degree
 sequences and extremal trees which attain the maximum value.
 \end{abstract}


\vskip 0.5cm

\section{Introduction}
The Wiener index of a molecular graph, introduced by Wiener
\cite{wiener1947} in 1947,
 is one of the oldest and most widely used topological indices in
the quantitative structure property relationships. In the
mathematical
 literature,  the Wiener index seems to be the first studied  by
 Entringer et al. \cite{entringer1976}.  For more information and background,
 the readers may refer to  a  recent and very
 comprehensive survey \cite{dobrynin2001} and a book \cite{rouvray2002}
 which is dedicated to Harry Wiener  on  the Wiener index and the references therein.

  Through this paper, all graphs are finite, simple and undirected.
    Let $G= (V,~E)$ be a simple connected graph with
vertex set $V(G)=\{v_1,\cdots, v_n\}$ and edge set $E(G)$.  Denote
by $d_G(v_i)$ (or for short $ d(v_i)$) the {\it degree} of vertex
$v_i$. The {\it distance} between vertices $v_i$ and $v_j$ is the
minimum number of edges between $v_i$ and $v_j$ and denoted by
$d_G(v_i, v_j)$ (or for short $d(v_i,v_j)$). The {\it Wiener index}
of a connected graph $G$ is defined as
\begin{equation}\label{weiner-def}
W(G)=\sum_{\{v_i, v_j\}\subseteq V(G)}d(v_i, v_j).
\end{equation}

A {\it tree} is a connected and acyclic graph.  A {\it caterpillar}
is a tree in which a single path (called {\it Spine}) is incident to
(or contains) every edge. For
 other terminology and notions, we follow from \cite{bondy1976}.

 Entringer et al. \cite{entringer1976} proved  that
  the path $P_n$  and the star $K_{1, n-1}$  have
the maximum  and minimum Wiener indices, respectively, in the set
consisting of all trees of order $n$. Dankelmann
\cite{dankelmann1994} obtained the all extremal graphs in the set of
all connected graphs with given  the order and the matching number
which attained the maximum  Wiener value. Moreover, Fischermann et
al. \cite{fischermann2002} and Jelen et al. \cite{jelen2003}
independently determined all trees which have the minimum Wiener
indices  among all trees of order $n$ and maximum degree $\Delta$.
 A nonincreasing
sequence of nonnegative integers
 $\pi=(d_1,d_2,\cdots, d_{n})$ is called {\it graphic} if there
 exists a  simple graph having $\pi$ as its  vertex degree sequence.
  Hence it is natural to
consider the following problem.
\begin{problem}\label{problem}
Let $\pi=(d_1, \cdots, d_n)$ be graphic degree sequence and
 $${\mathcal{G}}_{\pi}=\{G: {\rm \ the\  degree \  sequence\  of} \ G\ {\rm is} \ \pi\}.$$
 Find the upper (lower) bounds  for the Wiener index  of
all graphs $G$ in ${\mathcal G}_{ \pi}$ and characterize all
extremal graphs which attain the upper (lower) bounds.
\end{problem}
Moreover, we call a graph {\it maximum (minimum) optimal} if it
maximizes (minimizes) the Wiener index in $\mathcal{G_{\pi}}$.
Recently, by the different techniques, Wang \cite{wang2008} and
Zhang et al.\cite{zhang2008} independently characterized the tree
that minimizes the Wiener index among trees of given degree
sequences. Moreover, they proved that the  minimum optimal trees for
a given tree degree sequence $\pi$ are unique.
 On the other hand,
Wang in \cite{wang2008} also "{\it proved}" the only maximum optimal
tree that maximizes the Wiener index among trees of given degree
sequences. The result can be stated as follows:
\begin{theorem}\cite{wang2008}\label{wang}
Given the degree sequence and the number of vertices, the greedy
caterpillar maximizes the Wiener index, where the greedy caterpillar
with degree sequence $(d_1,\cdots, d_n)$ ($ d_1\ge d_2\ge \cdots \ge
d_k\ge 2>d_{k+1}=1$) is formed by attaching pending edges to a path
$v_1, v_2, \cdots, v_k$ of length $k-1$ such that
$$d(v_1)\ge d(v_k)\ge d(v_2)\ge d(v_{k-1})\ge \cdots\ge
d(v_{\lceil\frac{k+1}{2}\rceil}). $$
\end{theorem}

Unfortunately, this result is not correct. For example:
\begin{example}\label{example}
Let $\pi=(13, 5, 5, 5, 4, 3, 1, \cdots, 1)$ be a degree sequence of
tree with $31$ vertices.   Let $T_1$ and $T_2$ be two trees with
degree sequences $\pi$ (see Fig.1).

 \setlength{\unitlength}{0.1in}
\begin{picture}(60,15)
\put(3,5){\circle{0.5}} \put(3.25,5){\line(1,0){10}}
\put(13.5,5){\circle{0.5}}\put(13.75,5){\line(1,0){10}}
\put(24,5){\circle{0.5}} \put(24.25,5){\line(1,0){10}}
\put(34.5,5){\circle{0.5}} \put(34.75,5){\line(1,0){10}}
 \put(45,5){\circle{0.5}}
\put(45.25,5){\line(1,0){10}}
 \put(55.5,5){\circle{0.5}}

\put(2.8,5.2){\line(-1,2){3.5}} \put(3.2,5.2){\line(1,2){3.5}}

\put(13.5,5.25){\line(0,1){7}}

\put(23.8,5.2){\line(-1,2){3.5}} \put(24.2,5.2){\line(1,2){3.5}}

\put(34.3,5.2){\line(-1,2){3.5}} \put(34.7,5.2){\line(1,2){3.5}}

\put(44.8,5.2){\line(-1,2){3.5}} \put(45.2,5.2){\line(1,2){3.5}}

\put(55.3,5.2){\line(-1,2){3.5}} \put(55.7,5.2){\line(1,2){3.5}}

\put(-0.8,12.5){\circle{0.5}} \put(6.87,12.5){\circle{0.5}}

\put(2.2,12){$\cdots$}\put(2.2,13.3){$12$} \put(2.5,3.3){$v_{1}$}

\put(13.5,12.5){\circle{0.5}} \put(13,3.3){$v_{2}$}

\put(20.2,12.5){\circle{0.5}} \put(27.8,12.5){\circle{0.5}}
\put(23.2,13.3){$2$} \put(23.2, 3.3){$v_{3}$}

\put(30.61,12.5){\circle{0.5}} \put(38.3,12.5){\circle{0.5}}
\put(33.7,12){$\cdots$} \put(34.8,13.3){$3$} \put(34, 3.3){$v_{4}$}

\put(41.3,12.5){\circle{0.5}} \put(48.8,12.5){\circle{0.5}}
\put(44,12){$\cdots$} \put(44.8,13.3){$3$}\put(44, 3.3){$v_{5}$}

\put(51.7,12.5){\circle{0.5}}\put(59.2,12.5){\circle{0.5}}
 \put(54.5,12){$\cdots$}\put(55,13.3){$4$}\put(55, 3.3){$v_{6}$}
\put(28, 1){$T_1$}
          \end{picture}

\setlength{\unitlength}{0.1in}
\begin{picture}(60,15)
\put(3,5){\circle{0.5}} \put(3.25,5){\line(1,0){10}}
\put(13.5,5){\circle{0.5}}\put(13.75,5){\line(1,0){10}}
\put(24,5){\circle{0.5}} \put(24.25,5){\line(1,0){10}}
\put(34.5,5){\circle{0.5}} \put(34.75,5){\line(1,0){10}}
 \put(45,5){\circle{0.5}}
\put(45.25,5){\line(1,0){10}}
 \put(55.5,5){\circle{0.5}}

\put(2.8,5.2){\line(-1,2){3.5}} \put(3.2,5.2){\line(1,2){3.5}}

\put(13.3,5.2){\line(-1,2){3.5}} \put(13.7,5.2){\line(1,2){3.5}}

\put(23.8,5.2){\line(-1,2){3.5}} \put(24.2,5.2){\line(1,2){3.5}}

\put(34.5,5.25){\line(0,1){7}}

\put(44.8,5.2){\line(-1,2){3.5}} \put(45.2,5.2){\line(1,2){3.5}}

\put(55.3,5.2){\line(-1,2){3.5}} \put(55.7,5.2){\line(1,2){3.5}}

\put(-0.8,12.5){\circle{0.5}} \put(6.87,12.5){\circle{0.5}}

\put(2.2,12){$\cdots$}\put(2.2,13.3){$12$} \put(2.5,3.3){$v_{1}$}

\put(9.7,12.5){\circle{0.5}} \put(17.2,12.5){\circle{0.5}}
\put(12.6,12){$\cdots$} \put(13.5,13.3){$3$}
 \put(13,3.3){$v_{2}$}

\put(20.2,12.5){\circle{0.5}} \put(27.8,12.5){\circle{0.5}}
\put(23.2,13.3){$2$} \put(23.2, 3.3){$v_{3}$}

\put(34.5,12.5){\circle{0.5}}
 \put(34, 3.3){$v_{4}$}

\put(41.3,12.5){\circle{0.5}} \put(48.8,12.5){\circle{0.5}}
\put(44,12){$\cdots$} \put(44.8,13.3){$3$}\put(44, 3.3){$v_{5}$}

\put(51.7,12.5){\circle{0.5}}\put(59.2,12.5){\circle{0.5}}
 \put(54.5,12){$\cdots$}\put(55,13.3){$4$}\put(55, 3.3){$v_{6}$}
\put(20, 1){$T_2$} \put(25,-1){\bf Figure 1 $T_1$ and $T_2$}
          \end{picture}

\end{example}
Clearly, $T_2$ is a greedy caterpillar  and $T_1$ is not a greedy
caterpillar. Moreover, they have the same degree sequences $\pi$. By
calculation, it is easy to see that
$$W(T_2)=9870< W(T_1)=9886.$$
Hence this example illustrates that Theorem~\ref{wang} in
\cite{wang2008} is not correct.

Motivated by Problem \ref{problem} and Example \ref{example},  we
try to investigate the extremal trees which attain the maximum
Wiener index among all trees with given degree sequences. The
problem seems to be difficult. Because we find that the extremal
tree depends on the values of components of degree sequences.
 The rest of the paper is organized as follows. In Section 2,
  we discuss some properties of the extremal tree with the maximum Wiener
  index and give an upper bound in terms  of degree sequences.
 In Section
3, the extremal  trees with the maximum  Wiener index among given
degree sequences $(d_1, \cdots, d_n)$, where $d_1\ge \cdots \ge
d_k\ge 2>d_{k+1}=1$ and $k\le 6$ are characterized. Moreover, the
extremal maximal trees are not unique.

\section{Properties of extremal trees with the maximum Wiener index }

Let $\mathcal{T_{\pi}}$ be the set of all trees with degree
sequences $\pi=(d_1, d_2, \cdots,  d_n)$ with $d_1\ge
d_2\ge\cdots\ge d_n$.  Shi in \cite{shi1993} proved that a maximum
optimal tree must be a caterpillar.

\begin{lemma}\cite{shi1993}\label{shi}
Let $T^*$ be a maximum optimal tree in $\mathcal{T_{\pi}}$. Then
$T^*$ is a caterpillar.
\end{lemma}
From Lemma~\ref{shi}, we only need to consider all caterpillars with
a degree sequence $\pi$. In order to study the structure of the
maximum optimal trees, we present a formula for Wiener index of any
caterpillar.
\begin{lemma}\label{formula}
Let $T$ be a caterpillar of order $n$ with the degree sequence
$\pi=(d(v_1), \cdots,$ $ d(v_k), d(v_{k+1}),\cdots,d(v_n))$(see
Figure 2).

 \setlength{\unitlength}{0.1in}
\begin{picture}(60,20)
\put(10,5){\circle{0.5}} \put(10.25,5){\line(1,0){10}}
\put(20.5,5){\circle{0.5}}\put(20.75,5){\line(1,0){5}}

\put(26.5,4.55){$\cdot$}\put(27.5,4.55){$\cdot$}\put(28.5,4.55){$\cdot$}
\put(30,5){\circle{0.5}}
\put(31.5,4.55){$\cdot$}\put(32.5,4.55){$\cdot$}\put(33.5,4.55){$\cdot$}

\put(35.25,5){\line(1,0){5}}
\put(40.5,5){\circle{0.5}} \put(40.75,5){\line(1,0){10}}
 \put(51,5){\circle{0.5}}
\put(9.8,5.2){\line(-1,2){3.5}} \put(10.2,5.2){\line(1,2){3.5}}
\put(20.3,5.2){\line(-1,2){3.5}} \put(20.7,5.2){\line(1,2){3.5}}
\put(29.8,5.2){\line(-1,2){3.5}} \put(30.2,5.2){\line(1,2){3.5}}
\put(40.3,5.2){\line(-1,2){3.5}} \put(40.7,5.2){\line(1,2){3.5}}
\put(50.8,5.2){\line(-1,2){3.5}} \put(51.2,5.2){\line(1,2){3.5}}
\put(6.2,12.5){\circle{0.5}} \put(13.8,12.5){\circle{0.5}}
\put(16.65,12.5){\circle{0.5}} \put(24.25,12.5){\circle{0.5}}
\put(26.23,12.5){\circle{0.5}} \put(33.7,12.5){\circle{0.5}}
\put(36.7,12.5){\circle{0.5}} \put(44.23,12.5){\circle{0.5}}
\put(47.2,12.5){\circle{0.5}} \put(54.7,12.5){\circle{0.5}}

\put(9,12){$\cdots$}\put(19,12){$\cdots$}\put(29,12){$\cdots$}
\put(39,12){$\cdots$}\put(50,12){$\cdots$}

\put(9,13.3){$y_1$}\put(18.2,13.3){$y_{2}-1$}
\put(28.5,13.3){$y_{i}-1$} \put(37.5,13.3){$y_{k-1}-1$}
\put(50,13.3){$y_{k}$}

 \put(9.5,3.3){$v_{1}$}
\put(20,3.3){$v_{2}$} \put(30, 3.3){$v_{i}$} \put(40,3.3){$v_{k-1}$}
\put(50.2,3.3){$v_{k}$}
 \put(25,1){\bf Figure 2}\put(34, 1){$T$}
          \end{picture}
 If $d(v_i)=y_i+1\ge 2$ for $i=1, \cdots, k$ and
$d(v_{k+1})=\cdots=d(v_n)=1$, then
\begin{equation}\label{weiner-f}
W(T)=(n-1)^2+F(y_1, \cdots, y_k),
\end{equation}
where
\begin{equation}\label{fx}
F(y_1, \cdots,
y_k)=\sum_{i=1}^{k-1}(\sum_{j=1}^iy_j)(\sum_{j=i+1}^ky_j).
\end{equation}
 \end{lemma}
\begin{proof}
It is well known \cite{hosoya1971} that the formula
(\ref{weiner-def}) is equal to
$$W(T)=\sum_{e}n_1(e)n_2(e),$$
where $e=(u, v)$ is an edge of $T$, and $n_1(e)$ (resp. $n_2(e)$) is
the number of vertices of the component of $T-e$ containing $u$
(resp. $v$). For $e_i=(v_i, v_{i+1})\in E(T),$ the numbers of
vertices of the two components of $T-e_i$ are
$\sum_{j=1}^id(v_j)-(i-1)$ and $\sum_{j=i+1}^kd(v_j)-(k-i-1)$ for
$i=1, \cdots, k-1,$ respectively.  Hence
\begin{eqnarray*}
W(T)&=&\sum_{e\in E(T)}n_1(e)n_2(e)\\
&=&\sum_{e {\rm {\  is \  pendent\ edge}}}n_1(e)n_2(e)+
\sum_{ e  {\rm \  is \ not \ pendent\  edge}}n_1(e)n_2(e)\\
&=&
(n-1)(n-k)+\sum_{i=1}^{k-1}(\sum_{j=1}^id(v_j)-(i-1))(\sum_{j=i+1}^kd(v_j)-(k-i-1))
\\
&=&(n-1)(n-k)+\sum_{i=1}^{k-1}(1+\sum_{j=1}^iy_j)(1+\sum_{j=i+1}^ky_j)\\
&=&(n-1)(n-k)+(k-1)(1+\sum_{j=1}^ky_j)+\sum_{i=1}^{k-1}(\sum_{j=1}^iy_j)(\sum_{j=i+1}^ky_j)\\
&=&(n-1)^2+F(y_1, \cdots, y_k),
\end{eqnarray*}
where last equality is due to
$\sum_{j=1}^ky_j=\sum_{j=1}^kd(v_j)-k=2(n-1)-(n-k)-k=n-2$. This
completes the proof.
\end{proof}
{\bf Remark}  In this sequel, the caterpillar $T$ in
Lemma~\ref{formula} is denoted by $T(y_1,\cdots, y_k)$. Then degree
sequence of $T(y_1,\cdots, y_k)$ is $(y_1+1, \cdots, y_k+1, 1,
\cdots, 1)$. The following theorem give a characterization of a
maximum optimal tree.
\begin{theorem}\label{optimal-equal}
Let $\pi=(d_1, \cdots, d_n)$ with $d_1\ge\cdots \ge d_k\ge 2\ge
d_{k+1}=\cdots=d_n=1$. Then $T$ is a maximum optimal tree in
${\mathcal{T}}_{\pi}$ if and only if  $T$ is a caterpillar $T(x_1,
\cdots, x_k)$ and $(x_1, \cdots, x_k)$ satisfies
\begin{equation}
F(x_1, \cdots , x_k)=\max\{F(y_1, \cdots,
y_k)=\sum_{i=1}^{k-1}(\sum_{j=1}^iy_j)(\sum_{j=i+1}^ky_j):\  y_1\ge
y_k \},
\end{equation}
where $(y_1, \cdots, y_k)$ is any permutation of $(d_1-1, \cdots,
d_k-1)$.
\end{theorem}
\begin{proof} Necessity. Since $T$ is a maximum optimal tree in  ${\mathcal{T}}_{\pi}$, by
Lemmas~\ref{shi}, $T$ must be a caterpillar and can be denoted by
$T(z_1, \cdots, z_k)$ with $(z_1, \cdots, z_k)$ is the permutation
of $(d_1-1, \cdots, d_k-1)$. Moreover, by Lemma~\ref{formula}, we
have
$$W(T(z_1, \cdots, z_k))=(n-1)^2+F(z_1, \cdots, z_k).$$
For any permutation $(y_1, \cdots, y_k)$ of $(d_1-1, \cdots, d_k-1)$
with $y_1\ge y_k$, there exists a caterpillar $T_1$ with the degree
sequence $\pi$ such that
$$W(T_1)=(n-1)^2+F(y_1, \cdots, y_k).$$
Because $T(z_1, \cdots, z_k)$ is a maximum optimal tree in
${\mathcal{T}}_{\pi}$, we have
$$F(y_1, \cdots, y_k)=W(T_1)-(n-1)^2\le W(T(z_1, \cdots, z_k))-(n-1)^2=F(z_1, \cdots,
z_k).$$

Sufficiency. If $T$ is a caterpillar $T(x_1, \cdots, x_k)$ and
$(x_1, \cdots, x_k)$ satisfies
\begin{equation}
F(x_1, \cdots , x_k)=\max\{F(y_1, \cdots,
y_k)=\sum_{i=1}^{k-1}(\sum_{j=1}^iy_j)(\sum_{j=i+1}^ky_j): y_1\ge
y_k \},
\end{equation}
where the maximum is taken over all permutations $(y_1, \cdots,
y_k)$
 of $(d_1-1, \cdots, d_k-1)$. Let $T_1$ be any tree with the degree
sequence $\pi$. By Lemma~\ref{shi}, there exists a caterpillar $T_2$
with the degree sequence $\pi$ such that $W(T_1)\le W(T_2)$.   Then
$T_2$ must be $T(y_1, \cdots, y_k)$, where $(y_1, \cdots, y_k)$ is
the permutation of $(d_1-1, \cdots, d_k-1)$. Hence
$$W(T_1)\le W(T_2)= (n-1)^2+F(y_1, \cdots, y_k)\le (n-1)^2+F(x_1, \cdots,
x_k)=W(T(x_1, \cdots, x_k)).$$ Therefore $T(x_1, \cdots, x_k)$ is a
maximum optimal tree. This completes the proof.
\end{proof}

 Now we can present an upper bound for the Wiener index of any tree
with given degree sequence $\pi$ in terms of degree sequences.
\begin{theorem}\label{upperbound}
Let $T$ be a tree with a given  degree sequence $\pi=(d_1,\cdots,
d_n)$, where $d_1\ge \cdots\ge d_k>d_{k+1}=\cdots=d_n=1$. Then
\begin{equation}
W(T)\le (n-1)^2+\frac{k(k-1)}{4}\sum_{i=1}^k(d_i-1)^2
\end{equation}
with equality if and only if $k=2$ and $d_1=d_2$.
\end{theorem}
\begin{proof}
 Let $T(x_1, \cdots, x_k)$ be a
caterpillar and
 $(x_1, \cdots, x_k)$ satisfy
\begin{equation}
F(x_1, \cdots , x_k)=\max\{F(y_1, \cdots,
y_k)=\sum_{i=1}^{k-1}(\sum_{j=1}^iy_j)(\sum_{j=i+1}^ky_j): y_1\ge
y_k \},
\end{equation}
where $(y_1, \cdots, y_k)$ is any permutation of $(d_1-1, \cdots,
d_k-1)$. By Theorem~\ref{optimal-equal}, $W(T)\le W(T(x_1, \cdots,
x_k))$. Clearly,
\begin{eqnarray*} F(x_1,
\cdots, x_k)=\sum_{i=1}^{k-1}(\sum_{j=1}^ix_j)(\sum_{j=i+1}^kx_j)
=\frac{1}{2}(x_1, \cdots, x_k)C(x_1,\cdots, x_k)^T,
\end{eqnarray*}
where $$ C=\left(\begin{array}{cccccc} 0  &1&
2 &\cdots & k-2& k-1\\
1 & 0& 1& \cdots & k-3& k-2\\
\cdots & \cdots & \cdots & \cdots & \cdots & \cdots \\
k-1& k-2 & k-3&\cdots & 1 &0
\end{array}\right).$$
By Perron-Frobenius theorem (for example, see \cite{horn1985}), the
largest eigenvalue $\lambda_1(C)$ of $C$ is at most
$\frac{k(k-1)}{2}$ with equality if and only if $k=2$. Hence by
Rayleigh quotient,
$$(x_1, \cdots, x_k)C(x_1,\cdots, x_k)^T\le
\lambda_1(C)\sum_{i=1}^kx_i^2$$ with equality if and only if $(x_1,
\cdots, x_k)^T$ is an eigenvector of $C$ corresponding to the
eigenvalue $\lambda_1(C)$. Therefore,
$$F(x_1,
\cdots, x_k)\le \frac{k(k-1)}{4}\sum_{i=1}^kx_i^2$$ with equality if
and only if $k=2$ and $x_1=x_2$. Hence
$$ W(T)\le (n-1)^2+\frac{k(k-1)}{4}\sum_{i=1}^k{x_i}^2\le
(n-1)^2+\frac{k(k-1)}{4}\sum_{i=1}^k(d_i-1)^2$$ with equality if and
only if $k=2$ and $d_1=d_2$, since $(d(v_1), \cdots, d(v_k))$ is a
permutation of $(d_1, \cdots, d_k)$. This completes the proof.
\end{proof}

\begin{lemma}\label{function}
Let $w_1\ge w_2\ge\cdots\ge w_k\ge 1$ be the positive integers with
$k\ge 5$. Let
$$F(z_1, \cdots, z_k)=\max\{F(y_1, \cdots, y_k)=\sum_{i=1}^{k-1}(\sum_{j=1}^iy_j)(\sum_{j=i+1}^ky_j): y_1\ge y_k\},$$
where $(y_1, \cdots, y_k)$ is any permutation of $(w_1, \cdots,
w_k)$. Then there exists a $2\le t\le k-2$ such that the following
holds:
\begin{equation}\label{z1-zt}
z_1+\cdots+ z_{t-2}\le z_{t+1}+\cdots+z_k \end{equation} and
\begin{equation}\label{zt-zk}
z_1+\cdots+z_{t-1}>z_{t+2}+\cdots+z_k.
\end{equation}
Further, if equations (\ref{z1-zt}) is strict, then
\begin{equation}\label{z1-zt-ztrict}
z_1\ge z_2\ge\cdots\ge  z_t, \quad \quad z_t\le z_{t+1}\le \cdots\le
z_k.
\end{equation}
If equations (\ref{z1-zt}) becomes equality, then
\begin{equation}\label{lemma-z1-zt-1}
z_1\ge z_2\ge\cdots\ge  z_t, \quad \quad z_t\le z_{t+1}\le \cdots\le
z_k
\end{equation}
or
\begin{equation}\label{lemma-z1-zt-2}
z_1\ge z_2\ge\cdots\ge  z_{t-1}, \quad \quad z_{t-1}\le z_{t}\le
\cdots\le z_k.
\end{equation}
\end{lemma}
\begin{proof}
Let $$f(p)=\sum_{i=1}^{p-2}z_i-\sum_{i=p+1}^kz_i,\ \ \  2\le p\le
k-2.$$ Clearly $f(2)<0$, $f(k-1)>0$ and
$$f(2)\le f(3)\le\cdots\le f(k-1).$$
Hence there exists a $ 2\le t\le k-2$ such that $f(t)\le 0$ and
$f(t+1)>0$. In other words, equations (\ref{z1-zt}) and
(\ref{zt-zk}) hold. By the definition of $F(z_1, \cdots, z_k),$ we
have for $1\le i\le k-1$,
\begin{eqnarray*}
0&\le&  F(z_1, \cdots,z_{i-1}, z_i, z_{i+1}, \cdots, z_k)-F(z_1,
\cdots, z_{i-1},  z_{i+1}, z_i, \cdots, z_k)\\
&=&(z_{i+1}-z_i)(\sum_{j=1}^{i-1}z_j-\sum_{j=i+2}^kz_j) .
\end{eqnarray*}
But for $1\le i\le t-2$, by (\ref{z1-zt}), we have
$\sum_{j=1}^{i-1}z_j<\sum_{j=i+2}^kz_j$. Hence $z_1\ge \cdots\ge
z_{t-1}$. On the other hand,  for $t\le i\le k-1$, by (\ref{zt-zk}),
we have $\sum_{j=1}^{i-1}z_j>\sum_{j=i+2}^kz_j$. Therefore $z_t\le
z_{t+1}\cdots\le z_k$.

 If  (\ref{z1-zt}) is strict, then   $(z_1+\cdots+z_{t-2})-(z_{t+1}+\cdots+z_k)<0$,
 which implies
$z_{t-1}\ge z_t$. So (\ref{z1-zt-ztrict}) holds.

If (\ref{z1-zt}) becomes equality, i.e.,
$z_1+\cdots+z_{t-2}=z_{t+1}+\cdots+z_k$, then it is easy to see that
(\ref{lemma-z1-zt-1})  or (\ref{lemma-z1-zt-2}) holds. This
completes the proof.
\end{proof}

\begin{corollary}\label{k=6fun}
Let $w_1\ge w_2\ge\cdots\ge w_6\ge 1$ be the positive integers. Let
$$F(z_1, \cdots, z_6)=\max\{F(y_1, \cdots, y_6)=
\sum_{i=1}^{5}(\sum_{j=1}^iy_j)(\sum_{j=i+1}^6y_j):\ \ y_1\ge
y_6\},$$ where $(y_1, \cdots, y_6)$ is any permutation of $(w_1,
\cdots, w_6)$. Then $(z_1, \cdots, z_6)$  is equal to one of  the
following five $(w_1, w_6, w_5, w_4, w_3, w_2)$,  $(w_1, w_5, w_6,
w_4, w_3, w_2)$, $(w_1, w_4, w_6, w_5, w_3, w_2)$, $(w_1, w_4, w_5,
w_6, w_3, w_2)$ and $(w_1, w_3, w_6, w_5, w_4, w_2)$.
\end{corollary}
\begin{proof}
By Lemma~\ref{function}, there are just three cases:

{\bf Case 1} $t=2$. Then by Lemma~\ref{function}, $z_1\ge z_2$ and
$z_2\le z_3\le z_4\le z_5\le z_6$. Hence $(z_1, \cdots, z_6)$ must
be $(w_1, w_6, w_5, w_4, w_3, w_2)$.

{\bf Case 2} $t=3$. Then $z_1\le z_4+ z_5+z_6$ and
$z_1+z_2>z_5+z_6.$ Moreover,  $z_1\ge z_2\ge z_3$ and $ z_3\le
z_4\le z_5\le z_6$; or $z_1\ge z_2$ and $z_2\le z_3\le z_4\le z_5\le
z_6$.  Therefore $(z_1, \cdots, z_6)$ must be  one of  $(w_1, w_6,
w_5, w_4, w_3, w_2)$, $(w_1, w_5, w_6, w_4, w_3, w_2)$, $(w_1, w_4,
w_6, w_5, w_3, w_2)$ and $(w_1, w_3, w_6, w_5, w_4, w_2)$.

{\bf Case 3} $t=4$. Then $z_1+z_2\le z_5+z_6$. Moreover, $z_1\ge
z_2\ge z_3$ and $ z_3\le z_4\le z_5\le z_6$; or $z_1\ge z_2\ge
z_3\ge z_4$ and $z_4\le z_5\le z_6$. Therefore, $(z_1, \cdots, z_6)$
must be  one of  $(w_1, w_4, w_6, w_5, w_3, w_2)$,   $(w_1, w_5,
w_6, w_4, w_3, w_2)$ and $(w_1, w_4, w_5, w_6, w_3, w_2)$. This
completes the proof.
\end{proof}

\begin{theorem}\label{char} Let $\pi=(d_1, \cdots, d_n)$ be a  tree
degree sequence with $d_1\ge d_2\ge\cdots\ge d_k\ge 2$,
$d_{k+1}=\cdots=d_n=1$ and $k\ge 5$. If a caterpillar $T(x_1,
\cdots, x_k)$ is a maximum optimal tree in ${\mathcal{T}}_{\pi}$
with $F(x_1, \cdots, x_k)$ in equation (\ref{weiner-f}). Then there
exists a $2\le t\le k-2$ such that either
$$\sum_{i=1}^{t-2}x_i\le\sum_{i=t+1}^kx_i,\quad
\sum_{i=1}^{t-1}x_i>\sum_{t+2}^kx_i,\quad x_1\ge x_2\ge\cdots\ge
x_{t-1}\ge x_t, \quad  x_{t}\le x_{t+1}\le \cdots\le x_k;$$  or
$$\sum_{i=1}^{t-2}x_i=\sum_{i=t+1}^kx_i,\quad
\sum_{i=1}^{t-1}x_i>\sum_{t+2}^kx_i, \quad x_1\ge x_2\ge\cdots\ge
x_{t-1}\ge x_t, \quad  x_{t}\le x_{t+1}\le \cdots\le x_k;$$
 or
$$\sum_{i=1}^{t-2}x_i=\sum_{i=t+1}^kx_i,\quad
\sum_{i=1}^{t-1}x_i>\sum_{t+2}^kx_i, \quad x_1\ge x_2\ge\cdots\ge
x_{t-1}, \quad  x_{t-1}\le x_{t}\le \cdots\le x_k.$$
\end{theorem}
\begin{proof} It follows from Theorem~\ref{optimal-equal} and
Lemma~\ref{function} that the assertion  holds.
\end{proof}

\section{The maximum optimal tree with  many leaves}

 In this section, for a given degree sequence $\pi=(d_1, \cdots,
 d_n)$ with at least $n-6$ leaves, we give the maximum optimal trees with
 the maximum Wiener index in ${\mathcal{T}}_{\pi}$. Moreover, the
 maximum optimal tree may be not unique.
 \begin{theorem}\label{k=2--4}
 Let $\pi=(d_1, \cdots,d_k, \cdots, d_n)$ be tree degree sequence with
 $n-k$ leaves for $2\le k\le 4.$  Then the maximum optimal tree in ${\mathcal T}_{ \pi}$ is
 the greedy  caterpillar.
 In other words,

 if $k=2$, then  $W(T)=(n-1)^2+(d_1-1)(d_2-1)$, for $T\in {
 \mathcal{T}}_{\pi}.$

 If $k=3$, then for any $T\in {\mathcal{T}}_{\pi},$
  $$W(T)\le (n-1)^2+(d_1-1)(d_2+d_3-2)+(d_1+d_2-2)(d_3-1)$$
 with equality if and only if $T$ is the caterpillar $T(d_1-1, d_3-1,
 d_2-1).$

 If $k=4,$ then for  any $T\in {\mathcal{T}}_{\pi},$
 $$W(T)\le (n-1)^2+(d_1-1)(d_2+d_3+d_4-3)+(d_1+d_2-2)(d_3+d_4-2)+(d_1+d_2+d_3-3)(d_4-1)$$
 with equality if and only if $T$ is the caterpillar $T(d_1-1, d_4-1, d_3-1,
 d_2-1)$.
   \end{theorem}

\begin{proof} If $k=2$, it is obvious. If $k=3$,   it is easy to see
that  $F(d_1-1, d_2-1, d_3-1)\le F(d_1-1, d_3-1, d_2-1). $ By
Theorem~\ref{optimal-equal}, the assertion holds.

If $k=4$, then by Theorem~\ref{optimal-equal},  let $T$ be a
caterpillar $T(x_1, x_2, x_3, x_4)$ and
 $$F(x_1, x_2, x_3,
x_4)=\max\{F(y_1, y_2, y_3, y_4): y_1\ge y_4\},$$
 where $(y_1, y_2,
y_3, y_4)$ is any permutation of $(d_1-1, d_2-1, d_3-1, d_4-1)$.
Because
$$F(x_1, x_2, x_3, x_4)-F(x_2, x_1, x_3, x_4)=(x_1-x_2)(x_3+x_4)\ge
0$$ and $$ F(x_1, x_2, x_3, x_4)-F(x_1, x_2, x_4,
x_3)=(x_4-x_3)(x_1+x_2)\ge 0,$$ we have $x_1\ge x_2$ and $x_4\ge
x_3$. So $(x_1, x_2, x_3, x_4)=(d_1-1, d_4-1, d_3-1, d_2-1)$. This
completes the proof.
\end{proof}

\begin{theorem}\label{k=5}
 Let $\pi=(d_1, \cdots,d_k, \cdots, d_n)$ be tree degree sequence with
 $n-5$ leaves.

(1).  If $d_1> d_2+d_3$, then  the maximum optimal tree in
${\mathcal T}_{ \pi}$ is the only
 caterpillar $T(d_1-1, d_5-1, d_4-1, d_3-1, d_2-1)$.

(2). If $d_1=d_2+d_3$, then  there are the exactly two maximum
optimal trees in ${\mathcal T}_{ \pi}$:  one tree is the
 caterpillar  $T(d_1-1, d_5-1, d_4-1, d_3-1, d_2-1)$;
  the other tree is the caterpillar $T(d_1-1, d_4-1, d_5-1,
  d_3-1, d_2-1)$.

 (3). If $d_1< d_2+d_3$, then  the maximum optimal tree  in ${\mathcal T}_{ \pi}$ is the only
 caterpillar  $T(d_1-1, d_4-1, d_5-1, d_3-1, d_2-1)$.
 \end{theorem}
\begin{proof} By Theorem\ref{optimal-equal}, let $T(x_1, x_2, x_3, x_4, x_5)$
be a maximum optimal tree in ${\mathcal{T}}_{\pi}$.  If
$d_1>d_2+d_3$, then by Theorem~\ref{char}, it is easy to see that
$t=2$, and $x_1\ge x_2$ and $x_2\le x_3\le x_4\le x_5$. Hence $(x_1,
x_2, x_3, x_4, x_5)=(d_1-1, d_5-1, d_4-1, d_3-1, d_2-1)$.

If $d_1<d_2+d_3$, then by Theorem~\ref{char}, it is easy to see that
 $x_1\ge x_2\ge x_3$ and $x_3\le x_4\le x_5$. Hence $(x_1, x_2, x_3,
x_4, x_5)=(d_1-1, d_4-1, d_5-1, d_3-1, d_2-1)$  or $(d_1-1, d_3-1,
d_5-1, d_4-1, d_2-1)$. But $W(T(d_1-1, d_4-1, d_5-1, d_3-1,
d_2-1))-W(T(d_1-1, d_3-1, d_5-1, d_4-1,
d_2-1))=2(d_1-d_2)(d_3-d_4)\ge 0$ with equality if and only if
$d_1=d_2$ or $d_3=d_4$. Hence the assertion (3) holds.

If $d_1=d_2+d_3$, then by Theorem ~\ref{char}, it is easy to see
that  either $x_1\ge x_2$ and $x_2\le x_3\le x_4\le x_5$; or $x_1\ge
x_2\ge x_3$ and $x_3\le x_4\le x_5$. Hence $(x_1, x_2, x_3, x_4,
x_5)=(d_1-1, d_5-1, d_4-1, d_3-1, d_2-1)$ or $(d_1-1, d_4-1, d_5-1,
d_3-1, d_2-1)$. Moreover, $F(d_1-1, d_5-1, d_4-1, d_3-1, d_2-1)=F
(d_1-1, d_4-1, d_5-1, d_3-1, d_2-1)$. Hence  (2) holds.
\end{proof}

\begin{lemma}\label{6diff}
Let $w_1\ge w_2\ge\cdots\ge w_6\ge 1$ be positive integers and
$$ F(y_1, \cdots,
y_k)=\sum_{i=1}^{k-1}(\sum_{j=1}^iy_j)(\sum_{j=i+1}^ky_j).$$ Then
\begin{equation}\label{k61}
F(w_1, w_6, w_5, w_4, w_3, w_2)-F(w_1, w_5, w_6, w_4, w_3,
w_2)=(w_1-w_2-w_3-w_4)(w_5-w_6),
\end{equation}
\begin{equation}\label{k62}
F(w_1, w_5, w_6, w_4, w_3, w_2)-F(w_1, w_4, w_6, w_5, w_3,
w_2)=2(w_1-w_2-w_3)(w_4-w_5),
\end{equation}
\begin{equation}\label{k63}
F(w_1, w_4, w_6, w_5, w_3, w_2)-F(w_1, w_4, w_5, w_6, w_3,
w_2)=(w_1+w_4-w_2-w_3)(w_5-w_6),
\end{equation}
\begin{equation}\label{k64}
F(w_1, w_4, w_5, w_6, w_3, w_2)-F(w_1, w_3, w_6, w_5, w_4,
w_2)=(3w_3-3w_4-w_5+w_6)(w_1-w_2).
\end{equation}
\end{lemma}
\begin{proof}
By a simple calculation, it is easy to see that the assertion holds.
\end{proof}

\begin{theorem}\label{k=6}
 Let $\pi=(d_1, \cdots,d_6, \cdots, d_n)$ be tree degree sequence with
 $n-6$ leaves, i.e., $d_1\ge\cdots\ge d_6\ge 2$ and $d_7=\cdots=d_n=1$.

 (1). If $d_1>d_2+d_3+d_4-2$, then
there is only  one maximum optimal tree $T(d_1-1, d_6-1, d_5-1,
d_4-1, d_3-1, d_2-1)$ in ${\mathcal T}_{ \pi}$.

(2). If $d_1=d_2+d_3+d_4-2$, then
 there are exactly two maximum optimal trees in ${\mathcal T}_{
 \pi}$: one maximum optimal tree is $T(d_1-1, d_6-1, d_5-1,
d_4-1, d_3-1, d_2-1)$; the other maximum optimal tree is $T(d_1-1,
d_5-1, d_6-1, d_4-1, d_3-1, d_2-1)$.

(3). $d_2+d_3-1<d_1<d_2+d_3+d_4-2$, then there is only  one maximum
optimal tree $T(d_1-1, d_5-1, d_6-1, d_4-1, d_3-1, d_2-1)$ in
${\mathcal T}_{ \pi}$.

(4). If $d_2+d_3-1=d_1$, then there are exactly two maximum optimal
trees in ${\mathcal T}_{
 \pi}$: one maximum optimal tree is $T(d_1-1, d_5-1, d_6-1,
d_4-1, d_3-1, d_2-1)$; the other maximum optimal tree is $T(d_1-1,
d_4-1, d_6-1, d_5-1, d_3-1, d_2-1)$.

(5). If  $\  \ \max\{d_2+d_3-d_4,\
d_2+\frac{1}{3}(d_5-d_6)\}<d_1<d_2+d_3-1$, then there is only  one
maximum optimal tree $T(d_1-1, d_4-1, d_6-1, d_5-1, d_3-1, d_2-1)$
in ${\mathcal T}_{ \pi}$.

(6). If $d_1=d_2+d_3-w_4> d_2+\frac{1}{3}(d_5-d_6)$, then there are
exactly two maximum optimal trees in ${\mathcal T}_{
 \pi}$: one maximum optimal tree is $T(d_1-1, d_4-1, d_6-1,
d_5-1, d_3-1, d_2-1)$; the other maximum optimal tree is $T(d_1-1,
d_4-1, d_5-1, d_6-1, d_3-1, d_2-1)$.

(7). If $d_1= d_2+\frac{1}{3}(d_5-d_6)>d_2+d_3-d_4$, then there are
exactly two maximum optimal trees in ${\mathcal T}_{
 \pi}$: one maximum optimal tree is $T(d_1-1, d_4-1, d_6-1,
d_5-1, d_3-1, d_2-1)$; the other maximum optimal tree is $T(d_1-1,
d_3-1, d_6-1, d_5-1, d_4-1, d_2-1)$.

(8). If $d_1=d_2+d_3-d_4= d_2+\frac{1}{3}(d_5-d_6)$, then there are
exactly three maximum optimal trees in ${\mathcal T}_{
 \pi}$: they are  $T(d_1-1, d_4-1, d_6-1,
d_5-1, d_3-1, d_2-1)$; $T(d_1-1, d_4-1, d_5-1, d_6-1, d_3-1, d_2-1)$
and $T(d_1-1, d_3-1, d_6-1, d_5-1, d_4-1, d_2-1)$.

 (9). If $d_2+\frac{1}{3}(d_5-d_6)\le d_1<d_2+d_3-d_4$, or $d_1\le d_2+\frac{1}{3}(d_5-d_6)<
 d_2+d_3-d_4$,
  then there is only  one
maximum optimal tree $T(d_1-1, d_4-1, d_5-1, d_6-1, d_3-1, d_2-1)$
in ${\mathcal T}_{ \pi}$.

(10). If $d_2+d_3-d_4\le d_1<d_2+\frac{1}{3}(d_5-d_6)$; or $d_1\le
d_2+d_3-d_4< d_2+\frac{1}{3}(d_5-d_6)$,
 then there is only  one
maximum optimal tree $T(d_1-1, d_3-1, d_6-1, d_5-1, d_4-1, d_2-1)$
in ${\mathcal T}_{ \pi}$.

(11). If $d_1< d_2+\frac{1}{3}(d_5-d_6)=
 d_2+d_3-d_4$,
then there are exactly two maximum optimal trees in ${\mathcal T}_{
 \pi}$: one maximum optimal tree is $T(d_1-1, d_3-1, d_6-1,
d_5-1, d_4-1, d_2-1)$; the other maximum optimal tree is $T(d_1-1,
d_4-1, d_5-1, d_6-1, d_3-1, d_2-1)$.
\end{theorem}
\begin{proof}
 The proof is referred to appendix since it is technique.
\end{proof}
{\bf Remark}. From Theorem~\ref{k=6}, we can see that the maximum
optimal trees depend on the values of all components of the tree
degree sequences and not unique, while the minimum optimal tree is
unique for a given tree degree sequence. Moreover, Theorem~\ref{k=6}
explains that it seems to be difficult for characterize all the
maximum optimal trees for a given tree degree sequence.



\newpage
\begin{center}
Appendix:  Proof of Theorem~\ref{k=6}

\end{center}

\begin{lemma}\label{k=6fun}

Let $w_1\ge w_2\ge\cdots\ge w_6\ge 1$ be positive integers. If
$$F(z_1, \cdots, z_6)=\max\{F(y_1, \cdots, y_6): y_1\ge y_6\},$$
where $(y_1, \cdots, y_6)$  is any permutation of $(w_1, \cdots,
w_6)$, then the following statement holds.

(1). If $w_1>w_2+w_3+w_4$, then $(z_1, z_2, z_3, z_4, z_5,
z_6)=(w_1, w_6, w_5, w_4, w_3,w_2)$.

(2). If $w_1=w_2+w_3+w_4$, then $(z_1, z_2, z_3, z_4, z_5,
z_6)=(w_1, w_6, w_5, w_4, w_3,w_2)$ or $(w_1, w_5, w_6, w_4,
w_3,w_2)$.

(3). $w_2+w_3<w_1<w_2+w_3+w_4$, then $(z_1, z_2, z_3, z_4, z_5,
z_6)=(w_1, w_5, w_6, w_4, w_3,w_2)$.

(4). If $w_2+w_3=w_1$, then $(z_1, z_2, z_3, z_4, z_5, z_6)=(w_1,
w_5, w_6, w_4, w_3,w_2)$ or $(w_1, w_4, w_6, w_5, w_3,w_2)$.

(5). If  $\  \ \max\{w_2+w_3-w_4,
w_2+\frac{1}{3}(w_5-w_6)\}<w_1<w_2+w_3$, then $(z_1, z_2, z_3, z_4,
z_5, z_6)=(w_1, w_4, w_6, w_5, w_3,w_2)$.

(6). If $w_1=w_2+w_3-w_4> w_2+\frac{1}{3}(w_5-w_6)$, then $(z_1,
z_2, z_3, z_4, z_5, z_6)=(w_1, w_4, w_6, w_5, w_3,w_2)$ or $(w_1,
w_4, w_5, w_6, w_3,w_2)$.

(7). If $w_1= w_2+\frac{1}{3}(w_5-w_6)>w_2+w_3-w_4$, then $(z_1,
z_2, z_3, z_4, z_5, z_6)=(w_1, w_4, w_6, w_5, w_3,w_2)$ or $(w_1,
w_3, w_6, w_5, w_4,w_2)$.

(8). If $w_1=w_2+w_3-w_4= w_2+\frac{1}{3}(w_5-w_6)$, then $(z_1,
z_2, z_3, z_4, z_5, z_6)=(w_1, w_4, w_6, w_5, w_3,w_2)$, or $(w_1,
w_4, w_5, w_6, w_3,w_2)$ or $(w_1, w_3, w_6, w_5, w_4, w_2)$.

 (9). If $w_2+\frac{1}{3}(w_5-w_6)\le w_1<w_2+w_3-w_4$, or $w_1\le w_2+\frac{1}{3}(w_5-w_6)<
 w_2+w_3-w_4$, then
$(z_1, z_2, z_3, z_4, z_5, z_6)=(w_1, w_4, w_5, w_6, w_3,w_2)$.

(10). If $w_2+w_3-w_4\le w_1<w_2+\frac{1}{3}(w_5-w_6)$; or $w_1\le
w_2+w_3-w_4< w_2+\frac{1}{3}(w_5-w_6)$,  then $(z_1, z_2, z_3, z_4,
z_5, z_6)=(w_1, w_3, w_6, w_5, w_4,w_2).$

(11). If $w_1< w_2+\frac{1}{3}(w_5-w_6)=
 w_2+w_3-w_4$, then $(z_1, z_2, z_3, z_4,
z_5, z_6)= (w_1, w_4, w_5, w_6, w_3,w_2)$ or $(w_1, w_3, w_6, w_5,
w_4,w_2).$
\end{lemma}

\begin{proof}
(1). $w_1>w_2+w_3+w_4$. By (\ref{z1-zt}) and (\ref{zt-zk}) in
Lemma~\ref{function}, we have $t=2$ and $(z_1, \cdots, z_6)=(w_1,
w_6, w_5, w_4, w_3, w_2)$.

(2). $w_1=w_2+w_3+w_4$.  By (\ref{z1-zt}) and (\ref{zt-zk}) in
Lemma~\ref{function},  we have $t=3$. By (\ref{lemma-z1-zt-1}) and
(\ref{lemma-z1-zt-2}). we consider the following two cases.  If
$z_1\ge z_2\ge z _3$ and $z_3\le z_4\le z_5\le z_6$, then by
corollary~\ref{k=6fun} and $w_1=w_2+w_3+w_4$, we have $(z_1, \cdots,
z_6)=(w_1, w_5, w_6, w_4, w_3, w_2)$. If $z_1\ge z_2$ and $z_2\le
z_3\le z_4\le z_5\le z_6$, then $(z_1, \cdots, z_6)=(w_1, w_6, w_5,
w_4, w_3, w_2)$. Hence (2) holds.

(3). $w_2+w_3<w_1<w_2+w_3+w_4$.  We consider the following four
cases:

{\bf Case 1:} $w_2+w_3+w_5< w_1<w_2+w_3+w_4$.   By (\ref{z1-zt}) and
(\ref{zt-zk}) in Lemma~\ref{function}, we have $t=3$ and $z_1\ge
z_2\ge z _3$ and $z_3\le z_4\le z_5\le z_6$. Hence by
Corollary~\ref{k=6fun},  $(z_1, \cdots, z_6)=(w_1, w_5, w_6, w_4,
w_3, w_2)$.

{\bf Case 2:} $w_2+w_3+w_5=w_1<w_2+w_3+w_4$.  Similarly, $(z_1,
\cdots, z_6)=(w_1, w_5, w_6, w_4, w_3, w_2)$.

{\bf Case 3:} $w_2+w_4+w_5<w_1<w_2+w_3+w_5$ and $w_1>w_2+w_3$.    By
(\ref{z1-zt}) and (\ref{zt-zk}) in Lemma~\ref{function}, we have
$t=3$.  Further $(z_1, \cdots, z_6)=(w_1, w_5, w_6, w_4, w_3, w_2)$
or $(z_1, \cdots, z_6)=(w_1, w_4, w_6, w_5, w_3, w_2)$. But by
Lemma~\ref{6diff}, we have
$$F(w_1, w_5, w_6, w_4, w_3, w_2)-F(w_1, w_4, w_6, w_5, w_3,
w_2)=2(w_1-w_2-w_3)(w_4-w_5).$$ Hence  $(z_1, \cdots, z_6)=(w_1,
w_5, w_6, w_4, w_3, w_2)$.

{\bf Case 4:} $w_2+w_3<w_1\le w_2+w_4+w_5$.  By (\ref{z1-zt}) and
(\ref{zt-zk}) in Lemma~\ref{function}, we have $t=3$. Further $(z_1,
\cdots, z_6)=(w_1, w_3, w_6, w_5, w_4, w_2)$, or $(w_1, w_4, w_6,
w_5, w_3, w_2)$, or $(w_1, w_5, w_6, w_4, w_3, w_2)$. But by
Lemma~\ref{6diff}, we have
$$F(w_1, w_3, w_6, w_5, w_4, w_2)-F(w_1, w_2, w_6, w_5, w_4,
w_3)=2(w_2-w_3)(2w_1-w_4+w_6)\ge 0,$$
$$F(w_1, w_4, w_6, w_5, w_3, w_2)-F(w_1, w_3, w_6, w_5, w_4,
w_2)=(w_3-w_4)(3w_1-3w_2-w_5+w_6)\ge 0$$
 and
$$F(w_1, w_5, w_6, w_4, w_3, w_2)-F(w_1, w_4, w_6, w_5, w_3,
w_2)=2(w_1-w_2-w_3)(w_4-w_5).$$ Hence
 $(z_1,
\cdots, z_6)=(w_1, w_5, w_6, w_4, w_3, w_2)$.

(4). $w_2+w_3=w_1$.  From the proof of (3), it is easy to see that
$(z_1, \cdots, z_6)=(w_1, w_5, w_6, w_4, w_3, w_2)$ or $(w_1, w_4,
w_6, w_5, w_3, w_2)$, because $F(w_1, w_5, w_6, w_4, w_3,
w_2)-F(w_1, w_4, w_6, w_5, w_3, w_2)=0$. Therefore (4) holds.

(5). $\max\{w_2+w_3-w_4, w_2+\frac{1}{3}(w_5-w_6)\}<w_1<w_2+w_3$. We
consider the four cases.

{\bf Case 1:} $w_1>w_2+w_4+w_5$ and $w_1>w_2+w_3-w_5$. By
(\ref{z1-zt}) and (\ref{zt-zk}) in Lemma~\ref{function}, we have
$z_1\ge z_2\ge z_3$ and $z_3\le z_4\le z_5\le z_6$. Then $(z_1,
\cdots, z_6)=(w_1, w_5, w_6, w_4, w_3, w_2)$ or $(w_1, w_4, w_6,
w_5, w_3, w_2)$. But
$$F(w_1, w_5, w_6, w_4, w_3, w_2)-F(w_1, w_4, w_6, w_5, w_3,
w_2)=2(w_4-w_5)(w_1-w_2-w_3)\le 0$$ with equality if and only if
$w_4=w_5$. Therefore  $(z_1, \cdots, z_6)=(w_1, w_4, w_6, w_5, w_3,
w_2)$.

{\bf Case 2:} $w_1>w_2+w_4+w_5$ and $w_1\le w_2+w_3-w_5$. By
(\ref{z1-zt}) and (\ref{zt-zk}) in Lemma~\ref{function}, we have
$z_1\ge z_2\ge z_3$ and $z_3\le z_4\le z_5\le z_6$. Then $(z_1,
\cdots, z_6)=(w_1, w_5, w_6, w_4, w_3, w_2)$ or $(w_1, w_4, w_6,
w_5, w_3, w_2)$. But by Lemma~\ref{6diff}, we have
$$F(w_1, w_5, w_6, w_4, w_3, w_2)-F(w_1, w_4, w_6, w_5, w_3,
w_2)=2(w_4-w_5)(w_1-w_2-w_3)\le 0$$ with equality if and only if
$w_4=w_5$. Therefore  $(z_1, \cdots, z_6)=(w_1, w_4, w_6, w_5, w_3,
w_2)$.

{\bf Case 3:} $w_1\le w_2+w_4+w_5$ and $w_1> w_2+w_3-w_5$. By
(\ref{z1-zt}) and (\ref{zt-zk}) in Lemma~\ref{function}, we have
$z_1\ge z_2\ge z_3$ and $z_3\le z_4\le z_5\le z_6$. Then $(z_1,
\cdots, z_6)=(w_1, w_5, w_6, w_4, w_3, w_2)$, or $(w_1, w_4, w_6,
w_5, w_3, w_2)$, or $(w_1, w_3, w_6, w_5, w_4, w_2)$. But by
Lemma~\ref{6diff}, we have
$$F(w_1, w_5, w_6, w_4, w_3, w_2)-F(w_1, w_4, w_6, w_5, w_3,
w_2)=2(w_4-w_5)(w_1-w_2-w_3)\le 0$$ with equality if and only if
$w_4=w_5$. Moreover,
$$F(w_1, w_4, w_6, w_5, w_3, w_2)-F(w_1, w_3, w_6, w_5, w_4,
w_2)=(w_3-w_4)(3w_1-3w_2-w_5+w_6)\ge 0$$ with equality if and only
if $w_3=w_4$.
 Therefore  $(z_1, \cdots, z_6)=(w_1, w_4, w_6,
w_5, w_3, w_2)$.

{\bf Case 4:} $w_1\le w_2+w_4+w_5$ and $w_1\le w_2+w_3-w_5$. By
(\ref{z1-zt}) and (\ref{zt-zk}) in Lemma~\ref{function}, we have
$z_1\ge z_2\ge z_3$ and $z_3\le z_4\le z_5\le z_6$. Then  $(z_1,
\cdots, z_6)=(w_1, w_4, w_6, w_5, w_3, w_2)$, or $(w_1, w_3, w_6,
w_5, w_4, w_2)$. But by Lemma~\ref{6diff}, we have
$$F(w_1, w_4, w_6, w_5, w_3, w_2)-F(w_1, w_3, w_6, w_5, w_4,
w_2)=(w_3-w_4)(3w_1-3w_2-w_5+w_6)\ge 0$$ with equality if and only
if $w_3=w_4$. Therefore  $(z_1, \cdots, z_6)=(w_1, w_4, w_6, w_5,
w_3, w_2)$.

(6). $w_1=w_2+w_3-w_4> w_2+\frac{1}{3}(w_5-w_6)$. By (\ref{z1-zt})
and (\ref{zt-zk}) in Lemma~\ref{function}, we have $z_1\ge z_2\ge
z_3$ and $z_3\le z_4\le z_5\le z_6$; or $z_1\ge z_2\ge z_3\ge z_4$
and $ z_4\le z_5\le z_6$. Then  $(z_1, \cdots, z_6)=(w_1, w_4, w_6,
w_5, w_3, w_2)$;  or $(w_1, w_3, w_6, w_5, w_4, w_2)$; or$(w_1, w_5,
w_6, w_4, w_3, w_2)$; or $(w_1, w_4, w_5, w_6, w_3, w_2)$. But
$$F(w_1, w_4, w_6, w_5, w_3, w_2)-F(w_1, w_3, w_6, w_5, w_4,
w_2)=(w_3-w_4)(3w_1-3w_2-w_5+w_6)\ge 0.$$
$$F(w_1, w_5, w_6, w_4, w_3, w_2)-F(w_1, w_4, w_6, w_5, w_3,
w_2)=2(w_1-w_2-w_3)(w_4-w_5)\le 0.$$
$$F(w_1, w_4, w_6, w_5, w_3, w_2)-F(w_1, w_4, w_5, w_6, w_3,
w_2)=(w_1+w_4-w_2-w_3)(w_5-w_6)=0.$$ Therefore  $(z_1, \cdots,
z_6)=(w_1, w_4, w_6, w_5, w_3, w_2)$ or $(w_1, w_4, w_5, w_6, w_3,
w_2)$.

(7) $w_1= w_2+\frac{1}{3}(w_5-w_6)>w_2+w_3-w_4$. By (\ref{z1-zt})
and (\ref{zt-zk}) in Lemma~\ref{function}, we have $z_1\ge z_2\ge
z_3$ and $z_3\le z_4\le z_5\le z_6$. Then
 $(z_1, \cdots, z_6)=(w_1, w_4, w_6,
w_5, w_3, w_2)$;  or $(w_1, w_3, w_6, w_5, w_4, w_2)$; or $(w_1,
w_5, w_6, w_4, w_3, w_2)$. But by Lemma~\ref{6diff}, we have
$$F(w_1, w_4, w_6, w_5, w_3, w_2)-F(w_1, w_3, w_6, w_5, w_4,
w_2)=(w_3-w_4)(3w_1-3w_2-w_5+w_6)=0.$$
$$F(w_1, w_5, w_6, w_4, w_3, w_2)-F(w_1, w_4, w_6, w_5, w_3,
w_2)=2(w_1-w_2-w_3)(w_4-w_5)\le 0.$$  Hence $(z_1, \cdots,
z_6)=(w_1, w_4, w_6, w_5, w_3, w_2)$ or $(w_1, w_3, w_6, w_5, w_4,
w_2)$.

(8).  $w_1=w_2+w_3-w_4= w_2+\frac{1}{3}(w_5-w_6)$. It follows from
(6) and (7) that  (8) holds.

(9).  Assume that $w_2+\frac{1}{3}(w_5-w_6)\le w_1<w_2+w_3-w_4$. We
consider the following two cases:

{\bf Case 1:} $ w_1>w_2+w_4+w_5$. By (\ref{z1-zt}) and (\ref{zt-zk})
in Lemma~\ref{function}, we have $z_1\ge z_2\ge z_3$ and $z_3\le
z_4\le z_5\le z_6$; or  $z_1\ge z_2\ge z_3\ge z_4$ and $z_4\le
z_5\le z_6$. Hence
  $(z_1, \cdots, z_6)=(w_1, w_4, w_6,
w_5, w_3, w_2)$; $(w_1, w_4, w_5, w_6, w_3, w_2)$ or $(w_1, w_5,
w_6, w_4, w_3, w_2)$. But by Lemma~\ref{6diff}, we have
$$F(w_1, w_5, w_6, w_4, w_3, w_2)-F(w_1, w_4, w_6, w_5, w_3,
w_2)=2(w_1-w_2-w_3)(w_4-w_5)\le 0$$ with equality if and only if
$w_4=w_5$.
$$F(w_1, w_4, w_6, w_5, w_3, w_2)-F(w_1, w_4, w_5, w_6, w_3,
w_2)=2(w_1+w_4-w_2-w_3)(w_5-w_6)\le 0$$ with equality if and only if
$w_5=w_6$.  Therefore $(z_1, \cdots, z_6)=(w_1, w_4, w_5, w_6, w_3,
w_2)$.

 {\bf Case 2:} $ w_1\le w_2+w_4+w_5$. By (\ref{z1-zt}) and (\ref{zt-zk})
in Lemma~\ref{function}, we have $z_1\ge z_2\ge z_3$ and $z_3\le
z_4\le z_5\le z_6$; or $z_1\ge z_2\ge z_3\ge z_4$ and $z_4\le z_5\le
z_6$. Hence, $(z_1, \cdots, z_6)=(w_1, w_4, w_6, w_5, w_3, w_2)$; or
$(w_1, w_3, w_6, w_5, w_4, w_2)$; $(w_1, w_4, w_5, w_6, w_3, w_2)$
or $(w_1, w_5, w_6, w_4, w_3, w_2)$. But by Lemma~\ref{6diff}, we
have
$$F(w_1, w_4, w_6, w_5, w_3, w_2)-F(w_1, w_3, w_6, w_5, w_4,
w_2)=(3w_1-3w_2-w_5+w_6)(w_3-w_4)\le 0$$ with equality if and only
if $w_3=w_4$;
$$F(w_1, w_4, w_6, w_5, w_3, w_2)-F(w_1, w_5, w_6, w_4, w_3,
w_2)=2(-w_1+w_2+w_3)(w_4-w_5)\ge 0$$ and
$$F(w_1, w_4, w_6, w_5, w_3, w_2)-F(w_1, w_4, w_5, w_6, w_3,
w_2)=(w_1+w_4-w_2-w_3)(w_5-w_6)\ge 0$$ with equality if and only if
$w_5=w_6$. Therefore $(z_1, \cdots, z_6)=(w_1, w_4, w_5, w_6, w_3,
w_2)$.

Assume that $w_1\le w_2+\frac{1}{3}(w_5-w_6)<w_2+w_3-w_4$.    By
(\ref{z1-zt}) and (\ref{zt-zk}) in Lemma~\ref{function}, we have
$z_1\ge z_2\ge z_3$ and $z_3\le z_4\le z_5\le z_6$;  or $z_1\ge
z_2\ge z_3\le z_4$ and $z_4\le z_5\le z_6$. Hence, $(z_1, \cdots,
z_6)=(w_1, w_4, w_6, w_5, w_3, w_2)$;  or $(w_1, w_3, w_6, w_5, w_4,
w_2)$; or $(w_1, w_5, w_6, w_4, w_3, w_2)$;  or $(w_1, w_4, w_5,
w_6, w_3, w_2)$. But by Lemma~\ref{6diff}, we have
$$F(w_1, w_4, w_6, w_5, w_3, w_2)-F(w_1, w_3, w_6, w_5, w_4,
w_2)=(3w_1-3w_2-w_5+w_6)(w_3-w_4)\le 0$$ with equality if and only
if $w_3=w_4$;
$$F(w_1, w_5, w_6, w_4, w_3, w_2)-F(w_1, w_4, w_6, w_5,
w_3, w_2)=2(w_1-w_2-w_3)(w_4-w_5)\le 0;$$ and
$$F(w_1, w_4, w_5, w_6, w_3, w_2)-F(w_1, w_3, w_6, w_5,
w_4, w_2)=(3w_3-3w_4-w_5+w_6)(w_1-w_2)\ge 0$$ with equality if and
only if $w_1=w_2$. Therefore $(z_1, \cdots, z_6)=(w_1, w_4, w_5,
w_6, w_3, w_2)$.

(10). Assume that $w_2+w_3-w_4\le w_1<w_2+\frac{1}{3}(w_5-w_6)$.
  By
(\ref{z1-zt}) and (\ref{zt-zk}) in Lemma~\ref{function}, we have
$z_1\ge z_2\ge z_3$ and $z_3\le z_4\le z_5\le z_6$;  or $z_1\ge
z_2\ge z_3\le z_4$ and $z_4\le z_5\le z_6$. Hence, $(z_1, \cdots,
z_6)=(w_1, w_4, w_6, w_5, w_3, w_2)$;  or $(w_1, w_3, w_6, w_5, w_4,
w_2)$; or $(w_1, w_5, w_6, w_4, w_3, w_2)$;  or $(w_1, w_4, w_5,
w_6, w_3, w_2)$. But by Lemma~\ref{6diff}, we have
$$F(w_1, w_4, w_6, w_5, w_3, w_2)-F(w_1, w_3, w_6, w_5, w_4,
w_2)=(3w_1-3w_2-w_5+w_6)(w_3-w_4)\le 0$$ with equality if and only
if $w_3=w_4$;
$$F(w_1, w_5, w_6, w_4, w_3, w_2)-F(w_1, w_4, w_6, w_5,
w_3, w_2)=2(w_1-w_2-w_3)(w_4-w_5)\le 0;$$ and
$$F(w_1, w_4, w_6, w_5, w_3, w_2)-F(w_1, w_4, w_5, w_6,
w_3, w_2)=(w_1+w_4-w_2-w_3)(w_5-w_6)\ge 0$$ with equality if and
only if $w_5=w_6$. Therefore $(z_1, \cdots, z_6)=(w_1, w_3, w_6,
w_5, w_4, w_2)$.

 Assume that $w_1\le w_2+w_3-w_4<w_2+\frac{1}{3}(w_5-w_6)$.
  By
(\ref{z1-zt}) and (\ref{zt-zk}) in Lemma~\ref{function}, we have
$z_1\ge z_2\ge z_3$ and $z_3\le z_4\le z_5\le z_6$;  or $z_1\ge
z_2\ge z_3\le z_4$ and $z_4\le z_5\le z_6$. Hence, $(z_1, \cdots,
z_6)=(w_1, w_4, w_6, w_5, w_3, w_2)$;  or $(w_1, w_3, w_6, w_5, w_4,
w_2)$; or $(w_1, w_5, w_6, w_4, w_3, w_2)$;  or $(w_1, w_4, w_5,
w_6, w_3, w_2)$. But by Lemma~\ref{6diff}, we have
$$F(w_1, w_4, w_6, w_5, w_3, w_2)-F(w_1, w_3, w_6, w_5, w_4,
w_2)=(3w_1-3w_2-w_5+w_6)(w_3-w_4)\le 0$$ with equality if and only
if $w_3=w_4$;
$$F(w_1, w_5, w_6, w_4, w_3, w_2)-F(w_1, w_4, w_6, w_5,
w_3, w_2)=2(w_1-w_2-w_3)(w_4-w_5)\le 0;$$ and
$$F(w_1, w_4, w_5, w_6, w_3, w_2)-F(w_1, w_3, w_6, w_5,
w_4, w_2)=(3w_3-3w_4-w_5+w_6)(w_1-w_2)\le 0$$ with equality if and
only if $w_1=w_2$. Therefore $(z_1, \cdots, z_6)=(w_1, w_3, w_6,
w_5, w_4, w_2)$.

(11). $w_1< w_2+w_3-w_4=w_2+\frac{1}{3}(w_5-w_6)$. It follows from
$(9)$ and $(10)$ that (11) holds.
\end{proof}

\begin{proof} of Theorem~\ref{k=6}. It follows from Theorem~\ref{optimal-equal} and
Lemma~\ref{k=6fun} that the assertion holds.
\end{proof}
\end{document}